\documentclass[acmtog]{acmart}

\AtBeginDocument{%
  \providecommand\BibTeX{{%
    \normalfont B\kern-0.5em{\scshape i\kern-0.25em b}\kern-0.8em\TeX}}}

\setcopyright{acmlicensed}
\acmJournal{TOG}
\acmYear{2025} \acmVolume{44} \acmNumber{4} \acmArticle{} \acmMonth{8}\acmDOI{10.1145/3731179}

\usepackage{lineno}
\usepackage{layouts}
\usepackage{algorithm}
\usepackage{algpseudocode}
\usepackage[normalem]{ulem}
\usepackage{multirow}
\usepackage{makecell}
\usepackage{diagbox}
\usepackage{booktabs}
\usepackage{xcolor}
\usepackage{colortbl}

\algnewcommand\algorithmicinput{\textbf{Input:}}
\algnewcommand\algorithmicoutput{\textbf{Output:}}
\algnewcommand\algorithmicbl{\textbf{\# of blocks:}}
\algnewcommand\algorithmicth{\textbf{\# of threads / block:}}
\algnewcommand\algorithmicgl{\textbf{Global:}}
\algnewcommand\algorithmicinout{\textbf{In\textbackslash Out:}}

\algnewcommand\Input{\item[\algorithmicinput]}
\algnewcommand\Output{\item[\algorithmicoutput]}
\algnewcommand\Global{\item[\algorithmicgl]}
\algnewcommand\InOut{\item[\algorithmicinout]}
\algrenewcommand\algorithmicindent{1.0em}

\definecolor{dgreen}{rgb}{0.0, 0.7, 0.0}
\definecolor{dred}{rgb}{0.65, 0.16, 0.16}
\definecolor{dblue}{rgb}{0.16, 0.16, 0.65}
\definecolor{dgrey}{rgb}{0.16, 0.16, 0.16}
  
\newcommand{\prev}[1]{ \iffalse \textcolor{red}{ \textbf{Prev:} #1} \fi}

\newcommand{\ra}{$\rightarrow$}

\begin{document}

\title{Adaptive Algebraic Reuse of Reordering in Cholesky Factorizations with Dynamic Sparsity Patterns}





\author{Behrooz Zarebavani}
\affiliation{%
  \institution{University of Toronto}
  \country{Canada}}
\email{behrooz.zarebavani@gmail.com}

\author{Danny M. Kaufman}
\affiliation{%
  \institution{Adobe Research}
  \country{USA}}
\email{kaufman@adobe.com}

\author{David I.W. Levin}
\affiliation{%
  \institution{NVIDIA, University of Toronto}
  \country{Canada}}
\email{dlevin@nvidia.com}

  \author{Maryam Mehri Dehnavi}
\affiliation{%
  \institution{NVIDIA Research, University of Toronto},
  \country{USA}}
\email{mdehnavi@nvidia.com}

\renewcommand{\shortauthors}{Behrooz Zarebavani}


\begin{abstract}
We introduce Parth, a fill-reducing ordering method for sparse Cholesky solvers with dynamic sparsity patterns (e.g., in physics simulations with contact or geometry processing with local remeshing). Parth facilitates the selective reuse of fill-reducing orderings when sparsity patterns exhibit temporal coherence, avoiding full symbolic analysis by localizing the effect of dynamic sparsity changes on the ordering vector. We evaluated Parth on over 175,000 linear systems collected from both physics simulations and geometry processing applications, and show that for some of the most challenging physics simulations, it achieves up to 14x reordering runtime speedup, resulting in a 2x speedup in Cholesky solve time—even on top of well-optimized solvers such as Apple Accelerate and Intel MKL.

\end{abstract}




\begin{teaserfigure}
\includegraphics[width=\textwidth]{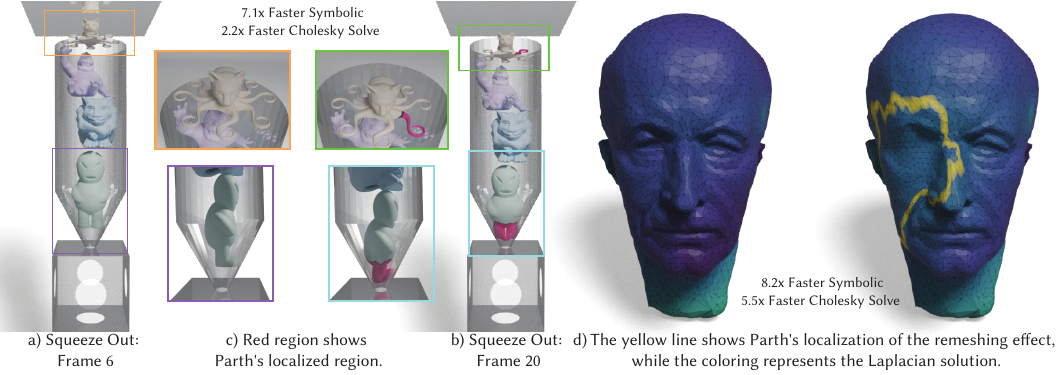}
 \caption{\textbf{Parth's adaptive reuse across changing sparsity patterns visualized for simulation and remeshing.} In (a) and (b), we demonstrate Parth's application in two frames of the large-deformation \textit{Squeeze Out} simulation benchmark from IPC~\cite{IPC}. Applying Parth to the simulation with the addition of three lines of code adaptively localizes these changes in the Hessian, with a resultant 7.1× speedup in symbolic analysis, leading to an overall 2.2× speedup for Cholesky solves. In Figure (c), we visualize Parth's localized change confinement on the simulation mesh, highlighted in red. Parth is also applied to a remeshing benchmark (d), in which we solve the Laplace–Beltrami operator on triangle meshes while remeshing. When the remeshing patch consists of 2\% of the mesh, Parth can reuse 92\% of the prior ordering computation, resulting in a 5.5× speedup in the Cholesky solve.}
\end{teaserfigure}
\maketitle


\section{Introduction}\label{sec:Introduction}

Linear solvers lie at the heart of many applications in graphics and scientific computing. Due to the structure induced by mesh-based computations, many common operations such as solving partial differential equations or minimizing a variational energy give rise to sparse systems of linear equations --  for the solution of which, it is either necessary or convenient to rely on Cholesky solvers. As a result, their accelerations are well-studied. However, the runtime of graphic applications often remains dominated by the cost of Cholesky solves, necessitating their further acceleration.

Efficient sparse Cholesky solvers are typically designed with static sparsity patterns in mind. The solution procedure itself is divided into two steps: "symbolic analysis," where the sparsity pattern of the matrix representing the system of linear equations is analyzed, and "numerical computation," which uses the symbolic analysis information to efficiently compute the solution. Often, a single symbolic analysis requires more computational resources than a single numerical computation (see Section~\ref{sec:Evaluation}). In cases of fixed sparsity, repeated solves can be accelerated by caching and reusing the results of a single symbolic computation. However, in applications with changing sparsity patterns, this optimization is unavailable.

In this work, we focus on 
optimizing the performance of Cholesky solvers in
applications with temporally coherent, local changes in sparsity pattern, such as those observed in contact simulation with elasticity~\cite{IPC} and geometric operations which involve local remeshing ~\cite{SurfaceMAp, OpenCloseSurface} in their computational pipeline. In such applications, the overhead of repeated
symbolic analyses becomes the bottleneck of linear solve costs. For example, in our evaluation of Incremental Potential Contact (IPC)~\cite{IPC}, we see that symbolic analysis accounts for up to 78\% of the total runtime of the Cholesky solver. Also, in our evaluation of the patch remeshing pipeline (see Section~\ref{fig:Evaluation:Remesher_SymBottleneck}), we observe that 82\% of the total runtime of the Cholesky solver is spent on symbolic analysis. By leveraging temporal coherence
in sparsity patterns, common in graphics workflows, we adaptively reuse prior symbolic analysis to accelerate sparse Cholesky solves.

Our main contribution, \textbf{Parth}, is an adaptive and general method that reuses symbolic analysis across repeated Cholesky solves when sparsity changes are localized and temporally coherent. Parth has two key objectives: (I) to enable adaptive reuse of symbolic analysis, and (II) to allow general-purpose integration with high-performance Cholesky solver libraries. These objectives ensure ease of use and portability while enabling speed-ups with state-of-the-art, reliable and performant Cholesky solvers.

Specifically, Parth accelerates the fill-reducing ordering step of symbolic analysis by reusing computation across calls to the symbolic analysis. Our evaluation of 175,000 linear systems from challenging physics simulations and geometric processing problems reveals that fill-reducing ordering is the primary bottleneck in symbolic analysis. This step computes a permutation to minimize fill-ins during the Cholesky solve, crucial for a fast solve~\cite{DirectSolverSurvey}. Focusing on the adaptive reuse of fill-reducing ordering allows us to integrate Parth into well-known performant Cholesky solvers as, while they differ in their underlying symbolic steps, they all employ fill-reducing ordering in their computational pipeline.

For example, by integrating Parth into Apple Accelerate~\cite{Apple}, and evaluating it on one of the challenging simulations in IPC (``Dolphin funnel''), we observe a 14x ordering runtime speedup without side effect on numerical performance. This speedup translates to an average of 2.2x speedup per solve, which also ends up accelerating the total solve time by 2x. Additionally, we demonstrate that by adding the three lines of code required to integrate Parth into the Cholesky solve computational pipeline of our IPC benchmark (see Section~\ref{sec:IPC:Benchmark}), our most challenging simulation (``Arma Roller'') achieves seven hours less computational time with only a 1.5× speedup in the total Cholesky solve runtime, without any side effects on numerical performance.

Parth achieves this speedup through a novel three-step process. First, it decomposes the graph dual of the input matrix using its hierarchical graph decomposition algorithm, which is suitable for decomposing fill-reducing orders. Next, it computes the fill-reducing order in each of the decomposed sub-graphs. Finally, when changes occur in the matrix representing a linear system—whether a change in the number of rows or columns or a change in the non-zero entries—Parth detects these changes and, using several novel algorithms, updates the decomposition to confine these changes and locally update the fill-reducing ordering vector. This process maintains an adaptive decomposition without the need for user tuning and allows for the reuse of computation in unchanged domains.

In summary, we present Parth, a method that reuses symbolic information in the presence of dynamic sparsity patterns with temporal coherence across calls to the Cholesky solver. To show its effectiveness, we evaluate Parth's performance across a wide range of linear systems with different behaviour in a change in the sparsity structure of the systems. We focus on three highest-performing Cholesky solvers—CHOLMOD~\citep{Cholmod}, the recently developed. 
Apple Accelerate sparse kernels~\cite{Apple}, and Intel MKL~\cite{MKLPardiso}. 
\footnote{We choose these by also comparing them with alternate available solvers, Sympiler~\cite{Sympiler}, Parsy~\cite{Parsy}, and Eigen. Please refer to Supplemental.}
We are open-sourcing the Parth codebase at \url{https://github.com/BehroozZare/Parth} to enable the community to readily integrate and benefit from our approach. Our extensive analysis, combined with public access to Parth, provides insights into how Parth improves performance across different solver frameworks and reveals specific performance characteristics of each state-of-the-art Cholesky solver. This transparency and accessibility will help practitioners make informed, application-specific decisions when selecting and tuning their preferred solver frameworks based on our comprehensive quantitative comparison.
\section{Related Work}\label{sec:RelatedWork}

Enhancing the performance of Cholesky factorizations and Sparse Triangular Solves (SpTRSV), remain a significant focus on computational mathematics and computer graphics due to their critical roles in simulation and optimization problems. 

\subsection{Exploiting Dense Computation}

Early optimizations center on exploiting dense computations within sparse factorizations. \citet{liu1990role} utilize the elimination tree to create supernode—groups of consecutive rows/columns sharing the same sparsity pattern. These supernodes are then factorized using dense BLAS\cite{BLAS3} kernels, enhancing computational efficiency by leveraging optimized dense linear algebra routines. Subsequent methods, such as CHOLMOD~\cite{Cholmod}, relax strict supernode constraints, allowing for better trade-offs by utilizing dense computation more effectively. These come at the expense of more redundant computation by increasing the size of the supernodes even when some rows/columns do not have a matching pattern.

\subsection{Parallelism and Scheduling Algorithms}

While the initial focus was primarily on finding dense computations within the chaos of sparse computations, later work focuses on further enabling parallelism across the computation of these dense blocks. This led to the introduction of advanced scheduling algorithms that balance load across parallel units, optimize memory usage through data reuse, and reduce synchronization overhead. MKL PARDISO~\cite{MKLPardiso} improves parallelism by efficiently distributing tasks across multiple cores and handling thread-level parallelism. Load-Balancing Coarsening (LBC)\cite{LBC} enhances memory usage by improving data locality and minimizing cache misses during computation on CPUs, as the parallelism of Cholesky factorization is limited by the sparsity pattern, and memory reuse can provide some compensation for the lack of parallelism. HDagg~\cite{zarebavani2022hdagg} provides a scheduler that balances the trade-off between memory reuse, parallelism and synchronization given hardware parallel capability. Additionally, the CHOLMOD GPU scheduler\cite{CholmodGPU} targets GPU hardware to achieve significant speedups in numerical computation, achieving up to a 2x speedup compared to CPU implementations. 


\subsection{Code Generation}

Recent approaches focus on providing optimized code by inspecting the sparsity pattern during symbolic analysis to improve numerical computation efficiency. Sympiler~\cite{Sympiler} analyzes sparsity patterns to generate specialized, pattern-specific code that optimizes memory access and computational efficiency. Similarly,~\citet{cheshmi2023runtime} automates the generation of high-performance code for sparse applications by fusing kernels used in linear solvers. However, this analysis comes with a high overhead of inspection, which further increases symbolic analysis overhead and must be redone if the sparsity pattern changes.

\subsection{Reuse of Numerical Factorization}

Efforts to reuse numerical factorization in linear solve computations focus on reducing the cost of recomputing factors for systems that exhibit numerical temporal coherence. In \citet{davis2005row}, a method is introduced to reuse factors through low-rank updates. \citet{tearing} extend this work to support low-rank updates when facing dynamic sparsity patterns. However, as they note, this approach does not perform well when multiple local changes occur in the sparsity pattern, since symbolic analysis cannot be reused under those circumstances. Furthermore, the method is not practical when the changes cannot be represented by low-rank vectors.

\begin{figure}[tp]
\centering
\includegraphics[width=\columnwidth]{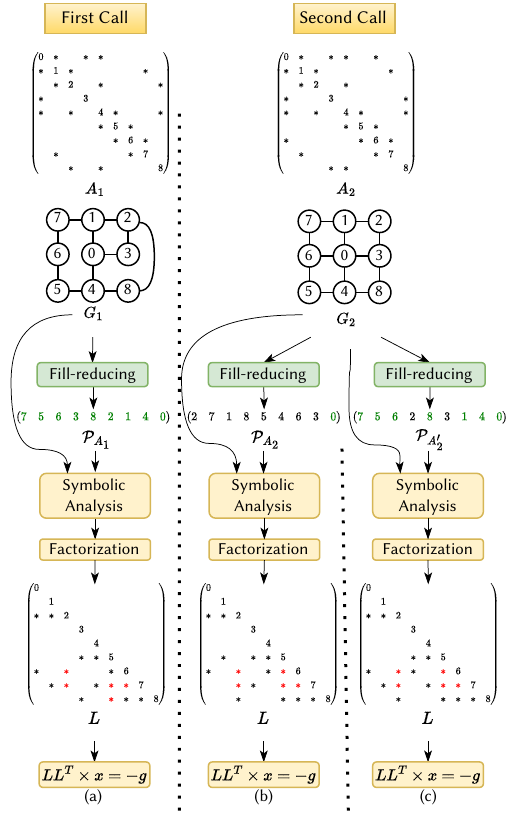}
\caption{\textbf{Cholesky Solver Pipeline Across Two Calls.} 
$A_1$ and $A_2$ are two consecutive linear systems with differing sparsity patterns, along with their corresponding graph duals $G_1$ and $G_2$—a reinterpretation of the system matrices as adjacency matrices. For example, removing the nonzero entries $(2,8)$ and $(8,2)$ from $A_1$ eliminates an edge in $G_1$. State-of-the-art fill-reducing ordering algorithms operate on such graphs. In the second call, applying fill-reducing ordering to $G_2$ can yield multiple high-quality permutations. Here, two possible orderings for $A_2$, $\mathcal{P}_{A_2}$ and $\mathcal{P}_{A'_2}$, are shown. Notably, $\mathcal{P}_{A'_2}$ shares many similarities (highlighted in green) with the previous ordering $\mathcal{P}_{A_1}$, illustrating the potential for reusing computations between calls while maintaining high-quality results.}
\label{fig:Background:Motive}
\end{figure}

Another body of work on reusing factorization computation is proposed in~\cite{FactorReuse1, FactorReuse2}, where a prior factor is reused across calls to the linear solver by only recomputing the affected supernodes. However, these methods often introduce approximation errors, limiting their applicability in contexts where the precision of Cholesky solvers is essential. Furthermore, they can not be applied to applications with dynamic sparsity patterns, as the elimination tree structure changes when the sparsity pattern changes across calls to Cholesky solvers, leading to the reconstruction of supernodes. NASOQ~\cite{cheshmi2020nasoq} provides a constraint-based QP solver that reuses factorization computation across calls to a constraint-based QP solver, as adding and removing these constraints requires a factorization. To handle changes in the sparsity pattern due to changes in constraints, it performs a full symbolic analysis with all possible constraints added and uses a subset of the symbolic analysis when the constraints are added or removed from the KKT matrix. However, this approach assumes prior knowledge of all potential constraints and can not handle new sparsity patterns not known in advance.

In contrast to existing methods, our work focuses on reusing symbolic analysis instead of numerical factorization to address the overhead caused by dynamic sparsity patterns in Cholesky solves. Specifically, we assume no prior knowledge of where changes may occur or how many changes can happen across each call to the Cholesky solver. By developing algorithms that identify and reuse unchanged portions of the symbolic analysis across iterations, we reduce the overhead associated with the symbolic phase. Our approach effectively handles multiple local changes while providing high-quality fill-reducing ordering.


\section{Background}\label{sec:ProblemStatement}

As solving the fill-reducing ordering problem is NP-hard~\cite{NPCompleteOrdering}, practical Cholesky solvers use heuristic methods such as METIS~\cite{METIS}, AMD~\cite{AMD}, and Scotch~\cite{Scotch} to generate high-quality orderings. These heuristics incorporate randomness, often producing multiple valid orderings for a single input. As an example, Figure~\ref{fig:Background:Motive} shows two systems of linear equations solved sequentially by a Cholesky solver with a gradual change between calls. For the second call, two high-quality ordering vectors are possible—one with high similarity to the previous call's order, indicating a solution that can be created by updating the prior ordering vector locally. This observation motivates Parth's approach.

Parth reuses computations from previous ordering module calls, narrowing the search space for high-quality orderings by leveraging the previously computed order. As a result, it generates a permutation vector with higher similarity to the prior vector than other approaches. 


Parth operates on the graph dual $G$ of the system of linear equations $A$ for a linear solve in the form of $Ax=b$, where $G$ is just a re-interpretation of $A$ as an adjacency matrix. This approach is identical to well-known ordering algorithms such as METIS and AMD which are graph-based. In this re-interpretation, each row/column corresponds to a node in graph $G$, and nonzero entries define the edges between nodes. This approach makes Parth general, as \( A \) could represent, for example, a Laplace-Beltrami operator or the Hessian of an energy function constructed for a Newton iteration solve. An example of graph duals $G_1$ and $G_2$ for matrices $A_1$ and $A_2$ is shown in Figure~\ref{fig:Background:Motive}. As this work focuses on mesh-based computations where the relationship between the graph and the underlying application, as described, is key.




\begin{figure}
\centering
\includegraphics{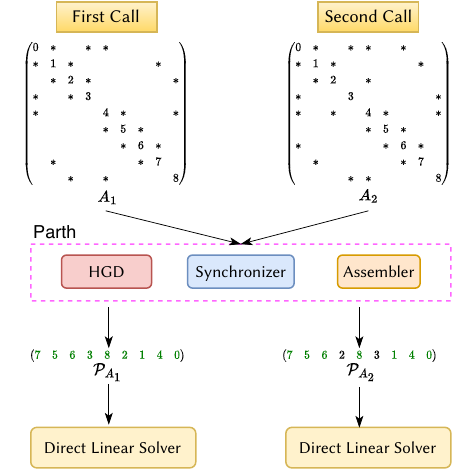}
\caption{\textbf{Parth's integration into high-performance Cholesky solvers for two calls with dynamic sparsity patterns.} The sparsity pattern shown here is the same as in Figure~\ref{fig:Background:Motive}. In the second call, Parth provides fill-reducing vectors with small changes by reusing computations from the first call and feeding them to the Cholesky linear solvers through their provided API. Note that current Cholesky linear solvers have no mechanism for reusing computations unless the sparsity pattern remains constant across calls. A detailed explanation of the HGD, Synchronizer, and Assembler modules, and how they interact, is provided in the rest of this section.}
\label{fig:Framework:Overview}
\end{figure}

\section{Parth Framework}\label{sec:Framework}

Parth replaces the fill-reducing module of  a sparse Cholesky solver and many implementations (Accelerate, MKL and CHOLMOD) provide such an interface. Parth consists of three underlying modules (Figure~\ref{fig:Framework:Overview}):  The \textit{\bf{Hierarchical Graph Decomposition}} (HGD) algorithm decomposes the dual graph $G$ of the system of linear equations $A$ into multiple smaller sub-graphs (fine-grain) that can be coarsened to form larger sub-graphs (coarse-grain). Furthermore, the decomposition enables the creation and combination of local permutation vectors per sub-graph. A \textit{\bf{Synchronizer}} module detects and integrates sparsity pattern changes into the created sub-graphs, thereby localizing the effects of these changes by confining them to specific coarse- or fine-grain sub-graphs. The \textit{\bf{Assembler}} unit then reuses the local permutation vectors in unchanged sub-graphs and updates those that are changed, finally assembling all this information into a single permutation vector, $\mathcal{P}_A$.

\subsection{Parth: Hierarchical Graph Decomposition}\label{sec:HGD}

The Hierarchical Graph Decomposition (HGD) algorithm is designed to enable the localization of changes within specific sub-graphs. The HGD module provides a decomposition that allows these changes to be confined to a single sub-graph or, if necessary, expands the decomposition to encapsulate broader changes. To ensure compatibility with fill-reducing ordering, HGD is constructed based on the nested dissection approach, as shown in Algorithm~\ref{alg:HGD}.


\begin{algorithm}[tp]
	\begin{algorithmic}[1]
        \Global $\mathcal{B}$
		\Input  $G_{sub}$, $l$, $i$, $max\_level$
		\If{ $l \neq max\_level$}
                \State  $g_l, g_r, g_s  \leftarrow computeMinSeparator(G_{sub})$
                \State $\mathcal{B}[i].nodes \leftarrow g_s.nodes$
                \Statex \textcolor{gray}{/*The left sub-graph in HGD procedure*/}
                \State $HGD(g_l, l + 1, 2 \times i + 1, max\_level)$
                \Statex \textcolor{gray}{/*The right sub-graph in HGD procedure*/}
                \State $HGD(g_r, l + 1, 2 \times i + 2, max\_level)$
            \Else
                \State $\mathcal{B}[i].nodes \leftarrow G_{sub}.nodes$
            \EndIf
	\end{algorithmic}
	\caption{Hierarchical Graph Decomposition ($HGD$)} 
	\label{alg:HGD}
\end{algorithm}

HGD (Algorithm~\ref{alg:HGD}) iteratively constructs a binary tree data structure, $\mathcal{B}$, stored as an array. The total number of nodes in the binary tree is computed using the maximum level parameter. A level is defined as the distance between the node and the root. Thus, HGD uses the $max\_level$ variable as the termination condition (line 1) for the recursion. Each call to the HGD algorithm results in a new binary tree node, represented as $\mathcal{B}[i]$. Each binary tree node corresponds to a sub-graph within the graph-dual of $A$. These subsets are determined by the recursive use of the function $computeMinSeparator$ (line 2). Note that this function divides the graph into three distinct sub-graph: a separator set, $g_s$, and two other sub-graphs, $g_l$ and $g_r$. The separator set acts as a minimal sub-graph that, when removed, dissects the remaining graph into two nearly equal parts, ensuring a balance in size between $g_l$ and $g_r$.


\begin{figure}
    \includegraphics{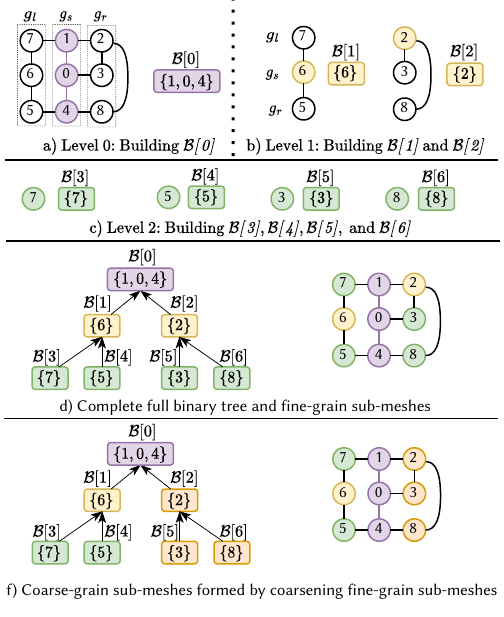}
    \caption{ {\bf Example of HGD evaluation on Figure~\ref{fig:Framework:Overview}:``First Call."} In this figure, the HGD algorithm operates across three levels as depicted in Figures~\ref{fig:FrameWork:HGD}(a-c). In Figure~\ref{fig:FrameWork:HGD}(a), the root node of the binary tree, $\mathcal{B}[0]$, is created, forming level 0 of the binary tree. The algorithm assigns the set $\{0,1,4\}$ to $g_s$, the set $\{5,6,7\}$ to $g_l$, and the set $\{2,3,8\}$ to $g_r$. Figures~\ref{fig:FrameWork:HGD}(b) and (c) illustrate the recursive expansion of the binary tree, generating all nodes at levels 1 and 2. Figure~\ref{fig:FrameWork:HGD}(d) displays the resulting fine-grain sub-graphs. Ultimately, Figure~\ref{fig:FrameWork:HGD}(f) reveals the coarse-grain sub-graph, which is derived from merging the separator $\mathcal{B}[2]$ with its ancestor nodes $\mathcal{B}[5]$ and $\mathcal{B}[6]$.}
\label{fig:FrameWork:HGD}
\end{figure}

Figure~\ref{fig:FrameWork:HGD} shows the HGD evaluation on \( G_1 \) in Figure~\ref{fig:Framework:Overview}:``First Call.'' Figures~\ref{fig:FrameWork:HGD}(a-c) present the HGD process at three levels. Initially, the root node of the binary tree is created, shown as \( \mathcal{B}[0] \). The decomposition then recursively advances to levels 1 and 2, resulting in the creation of 2 and 4 additional nodes within \( \mathcal{B} \). Each node in \( \mathcal{B} \) represents the smallest sub-graph in our decomposition. The final full binary tree is described in Figure~\ref{fig:FrameWork:HGD}(d). Note that each leaf of this binary tree represents a separate, approximately equal-sized sub-graph. The intermediate nodes in the binary tree are the separators, which also form part of the sub-graphs. This methodology effectively allows for coarsening sub-graphs by merging each sibling with its corresponding separator. Note that the coarsening can continue until all the sub-graphs are merged into the root, forming the whole \( G \). This demonstrates the hierarchical nature of the decomposition.

For a detailed explanation of why HGD is scalable and its computation is reusable, see the supplemental.


\begin{figure}
\centering
\includegraphics{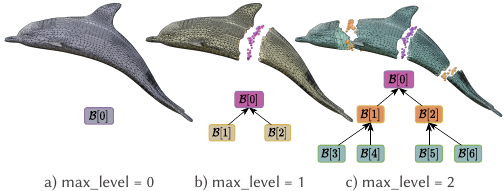}
\caption{{\bf HGD evaluation on "Dolphin" mesh for a Laplace-Beltrami operator.} Here, the graph \( G \) is identical to the mesh, as each node in the graph corresponds to a DOF in the mesh, and each edge in the graph arises from the connection of DOFs within an element. The figure shows the HGD evaluation for \( \text{max\_level} = \{0, 1, 2\} \). Note that in Figure~\ref{fig:Framework:HGDExampleRW}(b), the separator set is a small set of nodes compared to its left and right sub-graphs. Furthermore, observe how merging fine-grained sub-meshes \( \mathcal{B}[1] \), \( \mathcal{B}[3] \), and \( \mathcal{B}[4] \) in Figure~\ref{fig:Framework:HGDExampleRW}(c) results in the coarser sub-graph \( \mathcal{B}[1] \) in Figure~\ref{fig:Framework:HGDExampleRW}(b).}
\label{fig:Framework:HGDExampleRW}
\end{figure}

\subsection{Parth: Synchronizer}\label{sec:Synchronizer}

The Synchronizer module, in Algorithm~\ref{alg:Synchronizer}, synchronizes the information between $G$ and $\mathcal{B}$ across calls to the Cholesky solver by identifying changes in $G$ and encapsulating them within the set of sub-graphs in $\mathcal{B}$. If necessary, it re-decomposes the sub-graphs affected by these changes, as their nodes' connectivity may now differ. This allows Parth to preserve fill-reducing ordering information in sub-graphs that remain unchanged, enabling the reuse of this information. The Synchronizer achieves this through the five-step procedure.

 \begin{algorithm}[tp]
	\begin{algorithmic}[1]
            \Global $\mathcal{B}$
            \Input  $G$, $G_{new}$, $map$ $max\_level$, $Aggressive$
            \Output $C_\mathcal{B}$

            \State $setTrue(C_\mathcal{B})$
            \State $\mathcal{B} \leftarrow NodeChangeSynchronizer(map)$
            \State $E_G \leftarrow ComputeEdgeChanges(G, G_{new}, map)$
            \State $E_\mathcal{B} \leftarrow MapTo\mathcal{B}(E_G)$
            \Statex \textcolor{gray}{/* Fine- and coarse-grain dirty sub-graph*/}
            \State $D_F \leftarrow \{\}$, $D_C \leftarrow \{\}$ 
            \Statex \textcolor{gray}{/*Detect dirty sub-graphs in $\mathcal{B}$*/}
            \State $D_F, D_C \leftarrow DirtySubGraphDetection(E_\mathcal{B})$
            \Statex \textcolor{gray}{/*Filtering redundant work*/}
            \State FilterRedundantSubGraphs($D_F$, $D_S$)
            \Statex \textcolor{gray}{/*Mark the fine-grain sub-graphs for fill-reducing ordering*/}
            \State $C_\mathcal{B} \leftarrow MarkAndDecomposeSubGraphs(D_C, D_F)$
	\end{algorithmic}
	\caption{Synchronizer} 
	\label{alg:Synchronizer}
\end{algorithm}


\subsubsection{Step 1: Synchronizing the Added or Removed nodes}

The Synchronizer algorithm begins by analyzing whether the nodes in $G$ change, whether they are added or removed or their indices have changes (line 2). All this information can be computed from map input. The map array simply maps the index of nodes in the current call to the index of the same nodes in the previous call. For a detailed explanation of how $NodeChangeSynchronizer$ works, see supplemental. This function is used in our ``Remeshing'' benchmark. Based on our experience remeshers either create the map array explicitly in their underlying process or straightforwardly allow for its creation, as it is required for generating the faces and vertices metadata for (re)defining a mesh. Therefore, we hope this requirement is not too restrictive for other applications that could benefit from Parth.

\subsubsection{Step 2: Detecting Added or Removed Edges}

\begin{figure}
\includegraphics{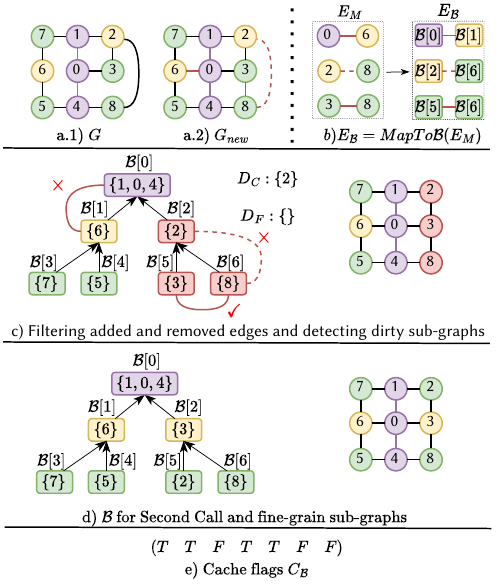}
    \caption{ {\bf Example of the Synchronizer procedure.} In Figure~\ref{fig:FrameWork:Synchronizer}(a,b), edge changes are detected and converted into edge changes across their corresponding fine-grain sub-graphs. In Figure~\ref{fig:FrameWork:Synchronizer}(c), the change that violates the separator condition (between $\mathcal{B}[5]$ and $\mathcal{B}[6]$) is considered, and other edge changes are disregarded. The lowest common descendant of $\mathcal{B}[5]$ and $\mathcal{B}[6]$ is $\mathcal{B}[2]$, where node 2 in $G_{new}$ no longer acts as a separator. The coarse-grain sub-graph encompassing the change is constructed by merging the sub-graphs in $\mathcal{B}[2]$, $\mathcal{B}[5]$, and $\mathcal{B}[6]$. In Figure~\ref{fig:FrameWork:Synchronizer}(d), the coarse-grain sub-graph represented by $\mathcal{B}[2]$ and its ancestors are re-decomposed. As a result, a new separator is chosen (node 3 in $G_{new}$). Figure~\ref{fig:FrameWork:Synchronizer}(e) displays the $C_\mathcal{B}$ array where 'T' (True) and 'F' (False) indicate which sub-graphs are intact and which have changed, respectively.}
\label{fig:FrameWork:Synchronizer}
\end{figure}

The Synchronizer now identifies added or removed edges from the graph by comparing the current graph ($G_{new}$) with the previous one ($G$) (line 3). Note that for an added DOF, all the edges are considered new, and a removed DOF's edges are not considered for computing the changed edges set $E_G$. After this step, edges in $E_G$ are mapped to changes in $\mathcal{B}$, indicated by $E_\mathcal{B}$. Figure~\ref{fig:FrameWork:Synchronizer}(a) illustrates the detection of changes between the graph duals shown in Figure~\ref{fig:Framework:Overview}. The difference between the two sub-graphs includes the addition of edges \textless0,6\textgreater and \textless3,8\textgreater and the deletion of edge \textless2,8\textgreater. Additionally, these modifications are mapped to the connectivity changes between fine-grain sub-graphs (line 4) represented by \textless$\mathcal{B}[0]$, $\mathcal{B}[1]$\textgreater, \textless$\mathcal{B}[2]$, $\mathcal{B}[6]$\textgreater, and \textless$\mathcal{B}[5]$, $\mathcal{B}[6]$\textgreater, as visualized in Figure~\ref{fig:FrameWork:Synchronizer}(b).

\begin{figure}
    \includegraphics{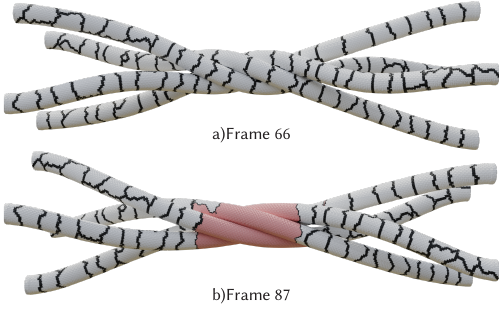}
    \caption{{\bf Synchronizer confinement of contacts in the ``Rods Twist'' simulation.} In this figure, multiple frames are skipped, and Parth is applied to the dual graph of the Hessian for the first Newton-solve iteration of the 66th and 87th frames. With four twisting rods, multiple contacts occur on each rod. To map sub-graphs to sub-graphs, every three consecutive nodes (3i, 3i+1, 3i+2) are mapped to a single DOF \(i\). Using \texttt{max\_level=7}, black lines represent separators, and white regions correspond to the leaves of \(\mathcal{B}\). The Synchronizer effectively coarsens small sub-graphes (in red) and confines 1612 changes, even with a step size of 21 frames.}
    \label{fig:Framework:SynchronizerRW}
\end{figure}

\subsubsection{Step 3: Detecting Dirty Sub-graphs}

The Synchronizer module utilizes $E_\mathcal{B}$ to identify sub-graphs affected by these edge changes. These affected sub-graphs, called "dirty sub-graphs," no longer possess valid information. There are two categories of dirty sub-graphs: $D_F$, which are fine-grain dirty sub-graphs, and $D_C$, which are coarse-grain dirty sub-graphs (line 5). The distinction is made because fine-grain sub-graphs do not require re-decomposition; only their fill-reducing ordering information needs to be updated. In contrast, coarse-grain sub-graphs may experience a violation of the separator's characteristics due to connectivity changes. That is, two sub-graphs that were supposed to be completely separated by a separator set are now connected via some nodes in the graph. Consequently, the Synchronizer must apply the HGD algorithm to re-decompose the coarse-grain regions and update the $\mathcal{B}$ information as needed. For a detailed explanation of $DirtySubGraphDetection$ algorithm (line 6), see supplemental.

Figure~\ref{fig:FrameWork:Synchronizer}(c) demonstrates the process of computing $D_F$ and $D_C$ using $E_\mathcal{B}$. The two changed edges \textless$\mathcal{B}[0]$, $\mathcal{B}[1]$\textgreater and \textless$\mathcal{B}[2]$, $\mathcal{B}[6]$\textgreater are disregarded because they do not disrupt the connection between the sub-graphs within the binary tree $\mathcal{B}$. For example, the change involving \textless$\mathcal{B}[0]$, $\mathcal{B}[1]$\textgreater is dismissed as the link between a separator and its adjacent left sub-graph does not violate the validity of the separator, given that $\mathcal{B}[0]$ separates $\mathcal{B}[1]$ from sub-graphs $\mathcal{B}[2]$, $\mathcal{B}[5]$, and $\mathcal{B}[6]$. Conversely, the change involving \textless$\mathcal{B}[5]$, $\mathcal{B}[6]$\textgreater violates the separator's role of $\mathcal{B}[2]$, since the sub-graph $\mathcal{B}[5]$ now connects to sub-graph $\mathcal{B}[6]$, resulting from the new connections between node 3 and node 8 in the graph. Consequently, the coarse sub-graph, annotated as $\mathcal{B}[2]$, comprising graph nodes $\{2,3,8\}$, is added to $D_C$ as it encapsulates the change \textless3,8\textgreater.

\subsubsection{Step 4: Filtering Sub-graphs}

During the formation of $D_F$ and $D_C$, each change is assessed independently from the others. As a result, some coarse-grained sub-graphs can encompass other fine- and coarse-grained sub-graphs. These smaller sub-graphs can subsequently be removed from $D_F$ and $D_C$, as they will be re-decomposed when the larger, encompassing sub-graph is re-decomposed. After this stage, the coarse sub-graphs requiring re-decomposition are fully identified. For a real-world example, refer to Figure~\ref{fig:Framework:SynchronizerRW}, where the contacts are restricted to a set of coarse sub-graphs (subsequently sub-graphs in the simulation mesh of IPC), colored in red, at the end of Step 4 (line 6 of Algorithm~\ref{alg:Synchronizer}).

\subsubsection{Step 5: Re-decomposing Sub-graphs}

Finally, after the creation and filtering of $D_F$ and $D_C$, the decomposition information in $\mathcal{B}$ must be updated accordingly. For fine-grain sub-graphs in $D_F$, we only need to mark them in the $C_\mathcal{B}$ array so that the Assembler module (Section~\ref{sec:Assembler}) updates the fill-reducing ordering information of these sub-graphs, as only the connectivity between the nodes has changed. For coarse-grain sub-graphs in $D_C$, the Synchronizer first re-decomposes the sub-graphs using the HGD algorithm (Section~\ref{sec:HGD}). The newly formed fine-grain sub-graphs are then marked in $C_\mathcal{B}$ for computation of fill-reducing ordering. For a detailed explanation of the $MarkAndDecomposeSubGraphs$, see supplemental.

As an illustration, on the right side of Figure~\ref{fig:FrameWork:Synchronizer}(c), the Synchronizer initially extracts the coarse sub-graph, colored in red (nodes ${2,3,8}$). Given that the size of $\mathcal{B}$ is constant, the Synchronizer is only required to substitute the invalid sub-graphs with valid ones. To achieve this, the process first identifies the sub-tree within $\mathcal{B}$ that contains all the invalid sub-graphs (illustrated in red on the left side of Figure~\ref{fig:FrameWork:Synchronizer}(c)). Next, it invokes the HGD algorithm on the sub-graph colored in red (Figure~\ref{fig:FrameWork:Synchronizer}(c)-right) to produce a valid decomposition, which is then visualized as a new sub-tree. This new sub-tree is used to replace the invalidated one, resulting in a completely valid binary tree $\mathcal{B}$, as depicted in Figure~\ref{fig:FrameWork:Synchronizer}(d). It is important to note that the sub-graph $\mathcal{B}[2]$ has been updated from the set ${2}$ to ${3}$, which now effectively isolates sub-graph $\mathcal{B}[5] = {2}$ from sub-graph $\mathcal{B}[6] = {8}$.

After updating the hierarchical graph decomposition, the new sub-graphs need the recomputation of fill-reducing ordering as the node structures they are representing have now changed. These dirty fine-grain sub-graphs are flagged within the $C_\mathcal{B}$ array. For instance, Figure~\ref{fig:FrameWork:Synchronizer}(e) indicates that the sub-graphs $\mathcal{B}[2]$, $\mathcal{B}[5]$, and $\mathcal{B}[6]$ are new and thus require updated fill-reducing information. As a result, both the $C_\mathcal{B}$ array and the binary tree $\mathcal{B}$ will be fed to the \textit{Assembler} module for new fill-reducing ordering information.

One issue is that a single edge connecting two sub-graphs with a common separator $\mathcal{B}[0]$ will result in the reuse of zero, as the coarse-grain sub-graph will encompass the whole graph. To alleviate this, we created a heuristic, named \emph{``Aggressive Reuse''}, explained in the supplemental. The heuristic moves one of the nodes that form the problematic edge to the corresponding separator. That is, the sub-graph forming that separator will expand if this happens, allowing Parth to achieve high reuse even in these scenarios.

\subsection{Parth: Assembler}\label{sec:Assembler}

The \textit{Assembler} module generates the permutation vector $\mathcal{P}_A$ by reusing computation across calls to Parth. Initially, the \textit{Assembler} computes the required permutation vectors for each sub-graph represented by $\mathcal{B}$. It then assembles these local permutation vectors to form $\mathcal{P}_A$ which applies to the whole matrix $A$. Due to changes in sparsity pattern, if \textit{Synchronizer} module marks some of the sub-graph "dirty" (refer to Section~\ref{sec:Synchronizer}), the \textit{Assembler} computes a new local permutation vector for these "dirty" sub-graphs. The values linked to these modified sub-graphs are then updated in $\mathcal{P}_A$, while the rest remain unchanged. This approach results in the computational reuse for all unchanged sub-graphs.

 \begin{algorithm}[tp]
	\begin{algorithmic}[1]
        \Global $\mathcal{B}$, $Order$
        \Input  $C_\mathcal{B}$, $DIM$
        \Output $\mathcal{P}_A$
        \State \textcolor{gray}{/*Initialization of $C_\mathcal{B}$*/}
        \If {$C_\mathcal{B}.empty()$} 
            \State $setFalse(C_\mathcal{B})$
            \State $Order \leftarrow PostOrdering(\mathcal{B})$
        \EndIf

        \State \textcolor{gray}{/*Computing offset of 
        $C_\mathcal{B}$*/}
        \State $Offset \leftarrow 0$
        \For{$i$ in $Order$}
                \State $\mathcal{B}[i].Offset \leftarrow Offset$
                \State $Offset \leftarrow Offset + |\mathcal{B}[i].nodes|$
        \EndFor
            
        \State \textcolor{gray}{/*Assembling fill-reducing ordering*/}
        \For{$i$ in $Order$}
            \If {$!C_\mathcal{B}[i]$} 
                \State $G_i \leftarrow getSubGraph(G, \mathcal{B}[i].nodes)$
                \State $\mathcal{B}[i].\mathcal{P}_l \leftarrow FillReducingOrdering(G_i)$
                \State $start \leftarrow \mathcal{B}[i].offset$
                \State $end \leftarrow start + |\mathcal{B}[i].nodes|$
                \State $\mathcal{P}_G[start:end] \leftarrow \mathcal{B}[i].\mathcal{P}_l$
            \EndIf
        \EndFor
        \Statex \textcolor{gray}{/*Convert mesh permutation into Hessian permutation*/}
        \For{$(j = 0;~j < |M.nodes|;~j = j + 1)$}   
            \For{$(d = 0;~d < DIM;~d = d + 1)$}   
                \State $\mathcal{P}_A[j * DIM + d] = \mathcal{P}_G[j] * DIM + d$
            \EndFor
        \EndFor
	\end{algorithmic}
	\caption{Parth: Assembler} 
	\label{alg:Assembler}
\end{algorithm}

\begin{figure}
\includegraphics{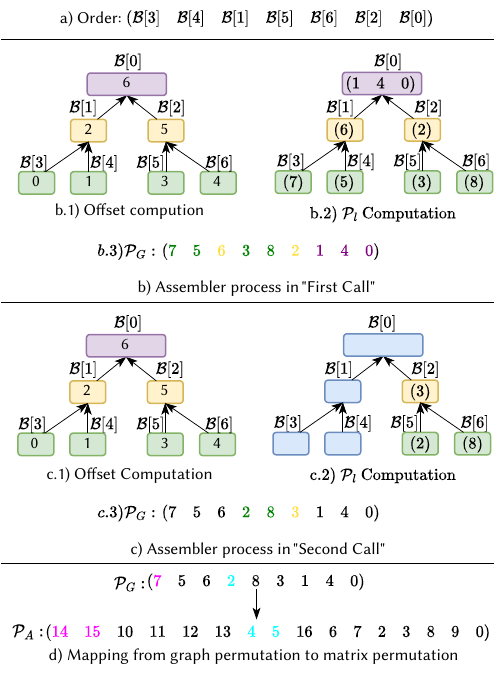}
    \caption{ {\bf Example of the Assembler procedure.} Considering two calls in Figure~\ref{fig:Background:Motive}, blue nodes in the $\mathcal{B}$ indicates computations that are skipped in ``Second Call'' due to reuse from ``First Call''. Figure~\ref{fig:FrameWork:Assembler}(a) displays the post-order traversal of $\mathcal{B}$, used in computing the position of each local permutation vector $\mathcal{P}_l$. During ``First Call'', the offsets, $\mathcal{P}_l$s, and the graph-based permutation vector $\mathcal{P}_G$ are computed. This initialization is presented in Figure~\ref{fig:FrameWork:Assembler}(b.1), (b.2), and (b.3). As detailed in Figure~\ref{fig:FrameWork:Synchronizer}, due to alterations in sub-graphs $\mathcal{B}[2]$, $\mathcal{B}[5]$, and $\mathcal{B}[6]$, the Assembler proceeds to update the offsets and $\mathcal{P}_l$ for these specific sub-graphs as shown in Figure~\ref{fig:FrameWork:Assembler}(c.1) and (c.2), and subsequently updates a portion of $\mathcal{P}_G$ depicted in Figure~\ref{fig:FrameWork:Assembler}(c.3). Finally, Figure~\ref{fig:FrameWork:Assembler}(d) illustrates the mapping of $\mathcal{P}_G$ to the $P_A$ for a 2D simulation.}
\label{fig:FrameWork:Assembler}
\end{figure}

The \textit{Assembler} procedure has three steps, as outlined in Algorithm~\ref{alg:Assembler}. In Step 1, Assembler determines the placement of the local permutation vectors $\mathcal{P}_l$ within $\mathcal{P}_A$. Step 2 involves the calculation of the necessary $\mathcal{P}_l$ vectors which, using the previously computed positions from Step 1, are inserted into the $\mathcal{P}_G$. Finally, in Step 3, the graph permutation vector $\mathcal{P}_G$ is converted into $\mathcal{P}_A$ based on $DIM$ value. All three steps are capable of reusing pre-computed data based on the array $C_\mathcal{B}$, which indicates the affected sub-graphes. Details of these three steps are further explained as follows:
 
\subsubsection{Step 1: Placement Computation} To use separators and their corresponding left and right sub-graphs as fill-reducing ordering information, we follow the nested-dissection approach. In this approach, the permutation vector is arranged in the computation so that the separator computation is placed after the left and right sub-graph computations (see the ``HGD Computation Reused'' section in supplemental for an example). Applying this order recursively is equal to placing each local permutation vector $\mathcal{P}_l$ related to each $\mathcal{B}$ node based on the post-order traversal of the tree. Consequently, the \textit{Assembler} module creates the entire graph permutation vector $\mathcal{P}_G$. The \textit{Assembler} executes this procedure in Lines 1-11 of Algorithm~\ref{alg:Assembler}.

Looking at Figure~\ref{fig:FrameWork:Assembler}, the post-order is computed once in Figure~\ref{fig:FrameWork:Assembler}(a), as the binary tree remains unchanged throughout the simulation. Subsequently, starting with $\mathcal{B}[3]$, the position of each $\mathcal{P}_l$ vector is identified, as illustrated in Figure~\ref{fig:FrameWork:Assembler}(b.3), based on the ``offset'' variable. For example, since sub-graphs $\mathcal{B}[3]$ and $\mathcal{B}[4]$ contain only a single node each, the aggregated node count resulting from visiting these sub-graph positions the $\mathcal{P}_l$ associated with sub-graph $\mathcal{B}[1]$ at the starting position $offset=2$ in $\mathcal{P}_G$. This procedure recurs in ``Second Call'' when the graph decomposition of three sub-graphs changes, leading to Figure~\ref{fig:FrameWork:Assembler}(c.3). Since computing the offset is not computationally intensive, we omit the reuse procedure of this step to simplify the explanation of offset computation.

\subsubsection{Step 2: $\mathcal{P}_l$ Computation and Mapping to $\mathcal{P}_G$} Once the placement of each $\mathcal{P}_l$ in $\mathcal{P}_G$ is determined, the \textit{Assembler} uses the array $C_\mathcal{B}$ to compute and update $\mathcal{P}_l$s for each sub-graph specified in $\mathcal{B}$ (lines 13-21). For instance, in Figure~\ref{fig:FrameWork:Assembler} (b.2 and c.2), the local permutation vectors $\mathcal{P}_l$ are shown for each sub-graph. In Figure~\ref{fig:FrameWork:Assembler}(c.2),`Second Call'', only three sub-graphs require new $\mathcal{P}_l$s, which allows Parth to reuse the $\mathcal{P}_l$ for the other four sub-graphs (coloured blue). By utilizing the ``offset'' variable calculated in Step 1, the Assembler constructs $\mathcal{P}_G$ by inserting the $\mathcal{P}_l$s into $\mathcal{P}_G$ (lines 13-21). Note that any ordering algorithm such as AMD~\cite{AMD}, Metis~\cite{Scotch}, and Morton Code can be used in $FillReducingOrdering$ (line 16), for $\mathcal{P}_l$ computation. Currently, Parth supports METIS and AMD for computing $\mathcal{P}_l$, but inserting new ones is straightforward.

\subsubsection{Step 3: Mapping $\mathcal{P}_G$ into $\mathcal{P}_A$} For 1D simulation ($DIM=1$), $\mathcal{P}_G$ is equal to $\mathcal{P}_A$. However, for different $DIM$, since Parth can compress the input $G$, or obtain the compressed form of $G$ by merging $DIM$ consecutive rows and columns in the form of $\{i * DIM + 1, \dots, i * DIM + DIM - 1\}$, the output $P_G$ is $DIM$ times smaller than $P_A$. Parth employs a straightforward mapping, as illustrated in lines 22-26, which results in forming $\mathcal{P}_A$. Note that this does not reduce quality, as these $DIM$ consecutive rows form a clique in $G$. Figure~\ref{fig:FrameWork:Assembler}(d) displays an example of mapping the graph simulation into $\mathcal{P}_A$ when $DIM=2$. Note that this mapping also utilizes reuse capability. However, to simplify the explanation, we've excluded the details of that implementation as these computations are fast compared to the $\mathcal{P}_l$ computation.


\section{Evaluation}\label{sec:Evaluation}

We focus on benchmarking Parth across a range of sparsity pattern variations, and  we evaluate Parth on both triangle and tetrahedral mesh geometries. To examine Parth's reuse performance when only non-zero entries change, we use Incremental Potential Contact (IPC)~\cite{IPC} simulations, referred to as the "IPC benchmark." To assess Parth in scenarios where the matrix dimensions change, we created a ``Remeshing benchmark''



For each benchmark, we first analyze the runtime bottleneck by determining how much time is spent on the numerical and symbolic phases of solves to illustrate the potential benefits of accelerating each part based on Amdahl's Law~\cite{Amdahl}. Furthermore, we demonstrate Parth's performance benefits when integrated into the solvers and how it reduces the symbolic analysis runtime. Finally, we demonstrate how Parth's performance impacts downstream numerical performance, highlighting the quality of Parth's fill-reducing ordering. In the following section, we summarize the hardware and software setup for our evaluations.

\subsection{Evaluation Setup}\label{sec:Evaluation:Setup}

Parth is implemented with C++. The separator computation uses METIS. For the local permutation vector, we allow the use of both AMD~\cite{AMD} and METIS. However, all of these libraries can be easily replaced as Parth's implementation does not rely on the underlying implementation of these ordering algorithms. For evaluation, we integrate Parth with, to our knowledge, the three most popular, highest-performing, robust Sparse Cholesky solver libraries: Intel MKL (MKL Pardiso LLT)~\cite{MKLPardiso}, SuiteSparse (CHOLMOD)~\cite{Cholmod}, and Apple Accelerate (Accelerate LLT)~\cite{Apple}. We do not include Eigen and Parsy~\cite{Parsy} Cholesky solvers, as they do not perform as well as the three high-performance libraries. See the ``Parsy and Eigen Evaluation'' section in the supplemental for more detail. We evaluate and compare the timing and accuracy between the Parth-augmented versions and the original versions of each of these linear solvers on Intel (20-core Xeon(R) Gold, 6248 CPU 2.5GHz, 28MB LLC cache, 202GB RAM) and Apple (12-core M2 Pro chip, 16GB RAM) platforms. For the Intel platform, we use MKL version 2023.4-912 and CHOLMOD version 7.6.0 with Ubuntu 22.04. For the Apple platform, we use the latest shipping Accelerate framework compiled with Xcode. We will release our evaluation code with this paper.

\begin{table}[tp]
    \caption{{\bf Summary of benchmark problem statistics.} We build a linear solve benchmark by precomputing and storing, in consecutive order, the sequential linear systems and geometries generated by the iterations of Newton solves for time-stepping six of the most challenging deformable-body benchmark problems from Li et al.\ \cite{IPC}. Here we summarize statistics for each simulation sequence with ID numbers for each (used throughout) on the left, \textit{$H_{rank}$} giving system dimension, and \textit{$H_{nnz}$} and \textit{$L_{nnz}$} respectively giving the average number of non-zeros per Hessian $H$ and corresponding $L$-factors (via METIS's symbolic analysis). The \#F represents a number of frames in each simulation while \textit{Changes} gives the percent of iterations per sequence where sparsity patterns change from a prior iteration.}
  \label{tab:Evaluation:SimProp}
  \rowcolors{2}{cyan!25}{white}
    \begin{tabular}{l c c c c} 
    \midrule
          Example & $H_{rank}$ & $H_{nnz}$, $L_{nnz}$ & \#F & Changes \\
     \hline
      (1) Dolphin Funnel & 24K & 500K, 11M & 800 & 99.9\%\\

     (2) Ball Mesh Roller   & 23K & 464K, 11M & 1000 & 96\%\\

     (3) Mat On Board  & 120K & 2.2M, 60M & 200 & 98.9\%\\

     (4) Rods Twist  & 160K & 3M, 54M & 1600 & 98.73\%\\

     (5) Squeeze out  & 135K & 2.8M, 72M & 1500 & 96\%\\

     (6) Arma Roller  & 201K & 4.8M, 344M & 400 & 99.9\% \\
    \bottomrule
    \end{tabular}

    \label{tab:Evaluation:SimProperties}
\end{table}

Each of the solver libraries mentioned above offers a range of options and default settings for fill-reducing ordering, significantly impacting solution quality and overall solve speeds. CHOLMOD's default configuration first applies AMD~\cite{AMD}. If a low-quality AMD ordering (based on measures of non-zeros and operation count) is generated, METIS~\cite{METIS} is then applied, and the better ordering of the two is adopted~\cite{Cholmod}. Across our benchmarks (see below), we observe that CHOLMOD almost always invokes METIS in 95.02\% of cases, with METIS being accepted as the final ordering. To improve CHOLMOD's overall performance, we configure it to directly apply METIS ordering in all subsequent experiments, eliminating this overhead from CHOLMOD's speeds. However, we emphasize that Parth also offers the option to use AMD when it is preferred for specific applications. MKL's LLT, by default, uses its own custom-optimized implementation of METIS, which we retain throughout our evaluation. Finally, Accelerate provides both AMD and METIS reordering options, with AMD as the default. As with CHOLMOD and as documented by Accelerate~\cite{Apple}, we observe a degradation in numerical performance speed when using AMD orderings on large-scale meshes compared to METIS in Accelerate's LLT. Thus, we likewise apply its METIS reordering for all benchmarks.

To ensure a fair comparison, we enhance all three libraries' solvers with additional logic to reuse symbolic analysis instead of recomputing it when the sparsity pattern remains unchanged across successive Cholesky solves in our benchmark. This eliminates runtime measurements for straightforward cases, ensuring that reported speedups reflect only when the sparsity patterns change.

\begin{table*}[tp]
    \caption{{\bf Breakdown of all Cholseky solves costs per library across all linear solves in our benchmark.} Here we summarize the timing (wall-clock seconds) and percent end-to-end Cholesky solver runtime per simulation breakdown per simulation sequence in our benchmark for all three state-of-the-art Cholesky solvers.}

    \centering
    \begin{tabular}{c|c c c c c c c}
    \midrule
     Tool & Step & Dolphin Funnel & Ball Mesh Roller & Mat On Board & Rods Twist & Squeeze out &  Arma Roller\\
     \hline
    \rowcolor{cyan!25}
     \cellcolor{white} & Symbolic time(s) & 3261 (77.8\%) & 2623 (70.7\%) & 1551 (76.1\%) & 8635 (75.5\%) & 9098 (78.1\%) & 35357 (53.5\%)\\
    \multirow{-2}{*}{ MKL} & Numeric time(s)& 932 (22.2\%) & 1085 (29.3\%) & 485 (23.9\%) & 2800 (24.5\%) & 2543 (21.9\%) & 30703 (46.5\%)     \\
    \midrule
    \midrule
        \rowcolor{cyan!25}
          \cellcolor{white}  & Symbolic time(s) & 1773 (71.52\%) & 1465 (62.97\%) & 985 (70.50\%) & 5363 (67.54\%) & 5540 (71.28\%) & 25331 (32.81\%)\\
         \multirow{-2}{*}{Accelerate} & Numeric time(s)& 706 (28.47\%) & 861 (37.03\%) & 412 (29.49\%) & 2578 (32.46\%) & 2231 (28.71\%) & 51871 (67.19\%)\\
    \midrule
    \midrule
            \rowcolor{cyan!25}
              \cellcolor{white}  & Symbolic time(s) & 2857 (40.5\%) & 2353 (34.3\%) & 1600 (49.03\%) & 8552 (48.4\%) & 9012 (48.7\%) & 35021 (37.5\%)\\
              \multirow{-2}{*}{CHOLMOD} & Numeric time(s) & 4200 (59.5\%) & 4511 (65.7\%) & 1663 (50.97\%) & 9098 (51.6\%) & 9483 (51.3\%) &  58282 (62.5\%)\\
    \midrule
    \end{tabular}

    \label{tab:Evaluation:Runtime}
\end{table*}

\subsection{IPC: Benchmark}\label{sec:IPC:Benchmark}

To analyze timings, bottlenecks and relative performance of these high-performance linear solvers in a consistent and fair side-by-side setting for IPC, we build a benchmark by precomputing and storing in consecutive order 143.5K sequential linear systems (Hessians, $A$, and gradients, $g$) which are generated by 5.5K Newton solves of challenging IPC volumetric FEM time-step problems~\cite{IPC}. We compute these systems by time-stepping six of the most challenging (over 96\% of all iterations in these simulations have changing sparsity due to contacts) deformable-body benchmark problems from~\citet{IPC}. Hessians in these systems range from 500K to 4.8M non-zero entries (with well over an order of magnitude increase in non-zeros for L-factors using state-of-the-art reordering with METIS~\cite{METIS}) depending on number of active contact stencils, model resolution and geometry; see Table~\ref{tab:Evaluation:SimProp}, and \citet{IPC} for simulation model statistics. We additionally store the initial conditions of each time step so that each Newton problem can also be analyzed consistently and independently across varying linear solvers. We use the IPC library~\cite{IPCCODE} to both generate this benchmark and to perform Newton solves for some of the analyses in the following sections.

Specifically, here we focus on providing insight into Parth's performance on volumetric, and triangular mesh (Mat On Board) simulation. Furthermore, this consecutive linear problem benchmark is critical for our analysis across high-accuracy linear solvers (within nonlinear-solve inner-loops) as convergence behaviour in Newton-type methods for stiff problems, as in the time-stepped elastodynamics simulations we consider here, are sensitive to minor changes in computed descent directions. These variations are generated by solvers due to rounding and parallelization and, in turn, produce differences in the numbers of linear system iterations per Newton solve (and so the entire simulation runs) when using different solvers. Note that, as we demonstrate in Supplemental (IPC: Numerical Performance Analysis), these variations do not generate iteration counts in favour of any of the Cholesky solvers. This benchmark then enables fair side-by-side evaluation of linear-solver methods within nonlinear solvers.

\subsection{IPC: Bottleneck Analysis}\label{sec:Evaluation:Bottleneck}

Table~\ref{tab:Evaluation:Runtime} summarizes the breakdown of total runtime costs per Cholesky solver library for the linear solves of each simulation sequence in the benchmark. Here, we see that for the two significantly faster solvers, MKL and Accelerate, symbolic analysis is clearly the primary bottleneck. For CHOLMOD, the story is a bit more nuanced. While CHOLMOD's symbolic computation runtimes are closely in line with MKL's, its significant slowdown compared to MKL is in its numeric computation phase, which ranges from two to four times slower than MKL. Here, symbolic analysis remains a significant cost (generally the same magnitude as numerical costs), and future optimizations of its numerical phase should bring its numerical costs down similarly to MKL's.

\begin{figure}[tp]
\centering
\includegraphics{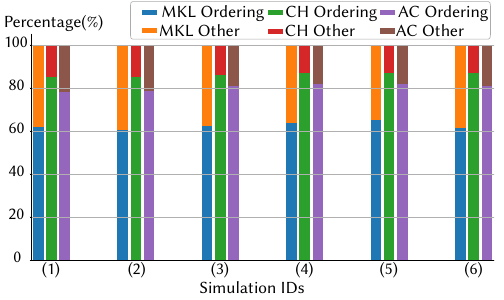}
\caption{{\bf Fill-reducing ordering is the bottleneck for all symbolic analyses.} Here we summarize a bottleneck analysis for the symbolic steps of all three Cholesky solvers across all six simulation sequences in our benchmark. For each simulation and corresponding tool, both the fill-reducing ordering time and all additional symbolic analysis time are recorded. Here, "CH", and "AC" respectively denote CHOLMOD and Accelerate LLT runtimes. "Other" categories summarize all other symbolic analysis costs per Cholesky solver with tasks that vary depending on the Cholesky solver method. 
}
\label{fig:Evaluation:SymBottleneck}
\end{figure}


This bottleneck, along with significant prior research and engineering focus in recent years on heavily optimizing numerical computation (see Section~\ref{sec:RelatedWork}), reiterates our research focus on improving symbolic analyses. Concretely, consider that MKL spends, on average, 71.8\% of its solve time on symbolic analysis. Best (an impossible zero-cost computation) improvements in numerical computation, e.g., the "Dolphin Funnel" sequence, would then be limited to a $1.28x$ speedup of linear solve costs.

If we then begin to look at the symbolic phase with finer granularity, we see that symbolic computation consists of multiple steps that vary with Cholesky solver implementation. However, \emph{all} sparse Cholesky solvers require and employ high-quality fill-reducing ordering methods in their symbolic phase. In turn, as we summarize in Figure~\ref{fig:Evaluation:SymBottleneck}, across all three Cholesky solvers, the fill-reducing ordering is by far the largest bottleneck in each solver's symbolic analysis steps. Here we see that across benchmark problems, fill-reducing ordering takes an average of $62.32\%$,  $81.05\%$, and $86.06\%$ of the total symbolic analysis time for MKL, Accelerate, and CHOLMOD respectively. At the same time, we observe that even the most opaque Cholesky solver libraries, which do not provide access to their symbolic analysis otherwise, offer APIs that can be used to replace default fill-reducing ordering implementations with customized methods. Thus, if we can show a significant boost in symbolic fill-reducing ordering, this offers the combined advantages of addressing the symbolic analysis bottleneck while staying modular to take advantage of the current, highly optimized, machine-specific numerical phases in existing high-performance solvers. 


\subsection{IPC: Parth Speedup}\label{sec:Evaluation:ReductionPerformance}

\begin{figure}[tp]
    \includegraphics[width=\columnwidth]{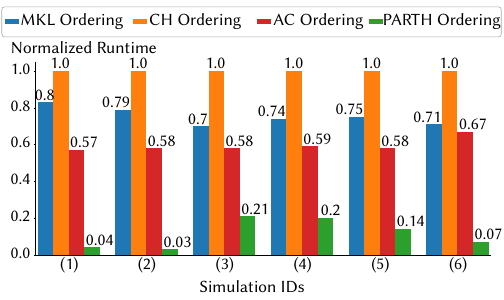}
    \caption{{\bf Parth Speedup.} the normalized runtime of Parth fill-reducing ordering compared to high-performance Cholesky solvers, namely, MKL, CHOLMOD (CH) and Accelerate (AC). Note that lower is better. The figure is normalized based on the slowest ordering algorithm. Across all 6 simulations, Parth fill-reducing ordering is faster than the fastest tool by 8.04x, 11.66x, 2.29x, 2.75x, 3.70x, and 9.79x speedup from simulation (1) to (6) respectively.}
    \Description{{\bf Parth Speedup.} the normalized runtime of Parth fill-reducing ordering compared to high-performance Cholesky solvers, namely, MKL, CHOLMOD (CH) and Accelerate (AC). Note that lower is better. The figure is normalized based on the slowest ordering algorithm. Across all 6 simulations, Parth fill-reducing ordering is faster than the fastest tool by 8.04x, 11.66x, 2.29x, 2.75x, 3.70x, and 9.79x speedup from simulation (1) to (6) respectively.}
\label{fig:Evaluation:ParthOrderingPerf}
\end{figure}

We first consider Parth's speedup compared to the fill-reduction timings of all three Cholesky solver libraries across our full benchmark, as Parth only accelerates this part of the linear solver pipeline. In Figure~\ref{fig:Evaluation:ParthOrderingPerf}, we summarize ordering runtimes normalized against the slowest solver. Here, we see that Parth outperforms all three Cholesky solvers across \emph{all} problem sequences in the benchmark, with speedups ranging from $2.8$X to well over an order of magnitude.

\begin{figure}
    \includegraphics[width=\columnwidth]{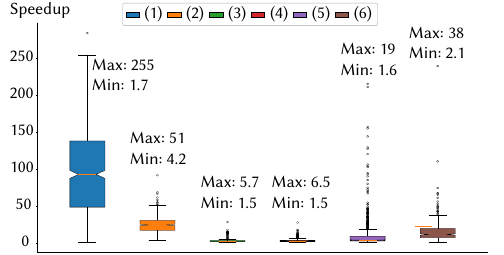}
    \caption{{\bf Parth speedup compared with best-available timing per solve.} The speedup achieved by Parth's fill-reducing ordering is compared to that of the best competitor among the three tools: Apple Accelerate, MKL, and CHOLMOD. Across all simulations, Parth's fill-reducing ordering consistently outperforms the state-of-the-art tool, resulting in speedups ranging from 1.5X to 255X (with outlying samples demonstrating further speedup in all sequences).}
\label{fig:Evaluation:OrderingWhisker}
\end{figure}

We next push this further and consider an ``optimal'' competing symbolic analysis which is enabled, without overhead, to pick the \emph{fastest} ordering among MKL, CHOLMOD and Accelerate (which can otherwise vary per solve for best speeds), for each Newton solve sequence in the benchmark. In Figure~\ref{fig:Evaluation:OrderingWhisker}, we compare Parth's speedup against this hypothetical best-speed analysis per Newton solve sequence. We find that even at the granularity of individual solves, Parth consistently remains significantly faster than the next best tool, with speedups per solve ranging from 1.5x to 255x across a wide range of solve sequences with both large and small sparsity changes. This is partly due to Parth's compression of the graph dual $G$, which can be coarsened when $DIM=3$. To be comprehensive, in supplemental, we provide a detailed analysis of Parth speedup on fill-reducing ordering when Parth does not employ compression. Finally, note that these speedups are reported only when the sparsity pattern changes, excluding the straightforward case of a constant sparsity pattern.

\subsection{IPC: Parth Ordering Quality}\label{sec:Evaluation:ReductionFill}

In the above analysis, we demonstrate that Parth efficiently computes permutation vectors with significant speedups in comparison to state-of-the-art symbolic methods across diverse sparsity patterns, Hessian sizes, and changing regions of sparsity updates. This demonstrates Parth's efficient reuse of computation across iterations within Newton solves, per time-step, and across sequential Newton solves, in time-stepped simulation sequences. With \emph{timing} settled we next analyze here the \emph{fill-reduction quality} of the orderings computed by Parth, and see that Parth delivers high-quality fill-reduction, comparable to the best sparsity generated among all three compared Cholesky solvers.

 \begin{figure}
    \includegraphics[width=\columnwidth]{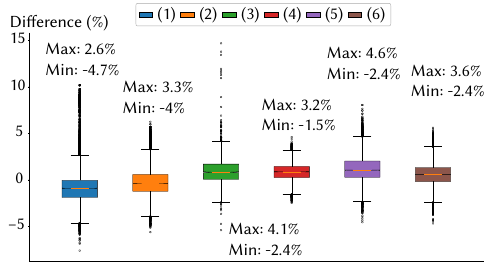}
    \caption{{\bf Fill-in reduction comparison across benchmark.} Here we we first identify the ``optimal'' fill-reducing ordering per solve, by choosing the sparsest factor generated across MKL, Accelerate, and CHOLMOD, for each. We then 
denote this smallest non-zero count among these three tools per linear system as $t_{best}$ and Parth's comparable non-zero count for its factor of the system as $t_{p}$. We then plot here and in Figure\ \ref{fig:Evaluation:LNNBallMeshRoller}
the difference $(t_{p} - t_{best}) / t_{best}$ quantifying the percent deviation between the non-zeros generated by the best high-performance Cholesky solvers and those by Parth per iteration.  Here the distribution of this measure across all simulation sequences shows that the median remains near zero, with minimum and maximum values generally (excluding outliers) falling well within $\pm 5\%$ of the best otherwise obtained solution.}
\label{fig:Evaluation:LNNZWhisker}
\end{figure}

In Figure~\ref{fig:Evaluation:LNNZWhisker} we compare Parth fill-reduction, with (as in the last section) an imagined, ``optimal''  competing symbolic analysis method across our benchmark. Specifically, we compare, per linear solve in the benchmark, against the \emph{best fill-reducing} ordering among MKL, CHOLMOD and Accelerate that generates the sparsest factor. Across all simulation sequences in the benchmark Parth's permutation quality remains consistent with the best sparse factors computed by the Cholesky solvers' fill-reduction—largely remaining well within a $\pm 5\%$ range, with medians close to zero. 

\begin{figure}
    \includegraphics[width=\columnwidth]{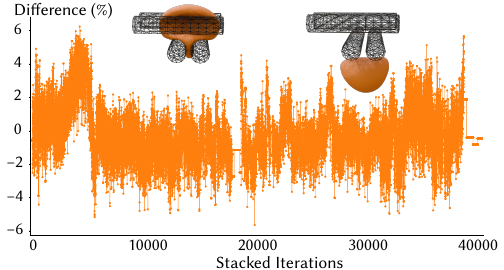}
    \caption{{\bf Detailed fill-in reduction analysis per solve.} The difference in the number of non-zeros (NNZ) in the factor generated by Parth is compared with that of the best results from MKL, Accelerate, and CHOLMOD. This comparison spans 40,000 stacked iterations from the Newton solves of time steps, with highly varying contact configurations and so large and rapid changes in sparsity patterns. The range of the difference is in $\pm6\%$, demonstrating Parth's generated sparse factors are comparable to the best available fill-reducing ordering while delivering large speedups in timings.}
\label{fig:Evaluation:LNNBallMeshRoller}
\end{figure}

Looking even closer at individual iterations, in Figure~\ref{fig:Evaluation:LNNBallMeshRoller} we plot the relative sparsity difference in non-zeros between the best-of solution and Parth across the first 40K successive linear solves in the high-contact and compression ``Ball Mesh Roller'' simulation sequence. Here we see that the total relative sparsity difference comparably remains in a range of $\pm 6\%$.

In summary, we find that across our benchmark, permutation vectors generated by Parth sometimes decrease, and sometimes slightly increase, sparsity (albeit both marginally) over the best provided by top competitor libraries. Recalling that heuristics applied in fill-reducing ordering, e.g. as in METIS~\cite{METIS}, lead to similar-scale minor variations in output when executed multiple times on the same problem, we see that Parth then provides comparable quality fill reduction at significant speedup. Finally, we also provide a comprehensive runtime analysis of numerical phase in supplemental which emphasizes the results shown in here.

 \begin{figure}
    \includegraphics[width=\columnwidth]{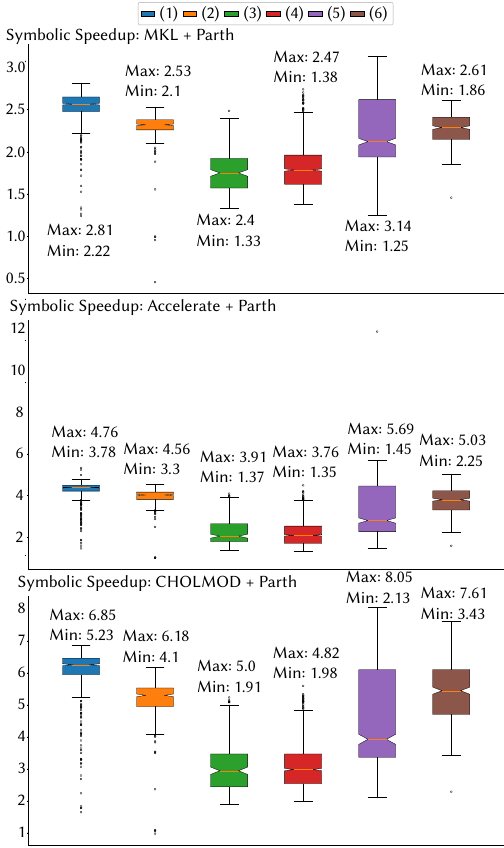}
    \caption{{\bf IPC: Parth symbolic performance impact across benchmark.} The three whisker plots indicate the performance impact of Parth on the symbolic step of  MKL, Accelerate and CHOLMOD. In all cases, except one iteration in ''Ball Mesh Roller", Parth improved the symbolic analysis performance by up to 8.05x speedup. }
\label{fig:Evaluation:SummaryWhisker_Sym}
\end{figure}

 \begin{figure}
    \includegraphics[width=\columnwidth]{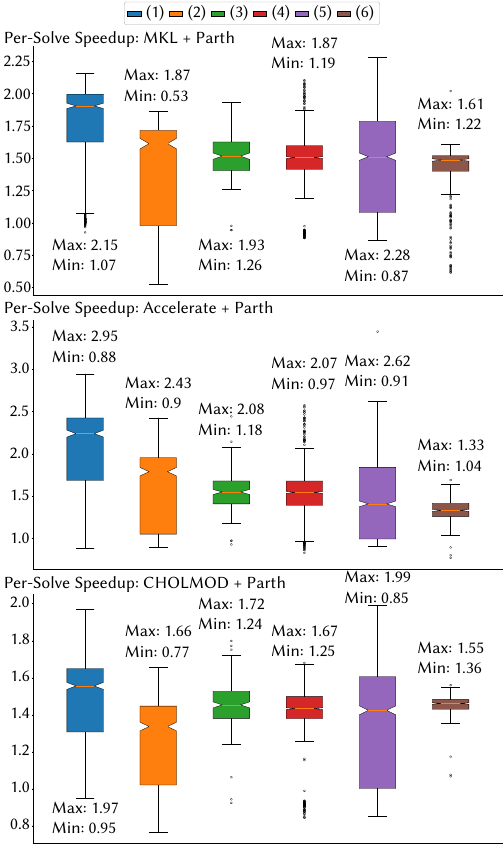}
    \caption{{\bf IPC: Parth's Per-Solve (Symbolic + Numeric) Performance Impact Across Benchmark.} As previously stated, MKL encountered numerical problems in the "Arma Roller"(6) simulation. Moreover, our analysis reveals that the minimum of 0.53x speedup recorded in "Roller Ball" (2) with MKL is due to a high number of repetitive iterations where the sparsity pattern remained unchanged. This implies a repeated use of a single permutation vector with slightly lower quality due to the absence of contact. Excluding these iterations, Parth achieves a 1.53 end-to-end speedup for this simulation.}
\label{fig:Evaluation:SummaryWhisker_perSolve}
\end{figure}

\subsection{IPC: Per Sparse Linear Solve Performance}\label{sec:Evaluation:ParthIntSolvePerformance}

Finally, we consider the performance of our three Parth-integrated Cholesky solvers on the symbolic and total costs of solves. Following the breakdown in Table~\ref{tab:Evaluation:Runtime}, we focus on the symbolic analysis phase per solver.
In Figure\ \ref{fig:Evaluation:SummaryWhisker_Sym}, across the benchmark, we see maximum speedups of 8.05x, 5.69x, and 3.14x for CHOLMOD, Accelerate, and MKL respectively, with corresponding median speedups of 4.9x, 3.5x, and 2.2x, and minimum speedups of 1.9x, 1.4x, and 1.2x.

In Figure \ref{fig:Evaluation:SummaryWhisker_perSolve}, we see that by just directly integrating the Parth module into Cholesky solver packages, we obtain speedups on full Cholesky solve costs of up to 2.95X, 2.28X, and 1.97X (1.5X, 1.54X, and 1.43X median) for Accelerate, MKL, and CHOLMOD respectively. In line with our earlier analysis, we also confirm that Parth consistently extracts the greatest performance speedups as it is integrated with successively more optimized and performant Cholesky solvers leveraging hardware—here in Accelerate and then with MKL. In Table~\ref{tab:Evaluation:ParthEffectedRuntime}, we report the \emph{total} Cholesky-solver runtime costs across all linear solves for each simulation sequence in the benchmark, and again observe consistent breakdowns.


\begin{table*}[tp]
    \caption{{\bf IPC: Breakdown of Cholesky solve costs per library across all Parth-integrated Cholesky solves in our benchmark.} Here we summarize the timing (wall-clock seconds) and present total linear solve runtime for default and Path-integrated versions per simulation sequence in our IPC benchmark for all three state-of-the-art Cholesky solvers corresponding to original costs in Table~\ref{tab:Evaluation:Runtime}. In here, the solve speedup from left (Dolphin Funnel) to right (Arma Roller) for 3 Cholesky solvers are as follow: MKL speedup=$\{1.79,1.68,1.49,1.56,.1.71,1.37\}$ and Accelerate speedup=$\{2.07,1.89,1.56,1.78,1.33\}$ and CHOLMOD speedup=$\{1.5,1.41,1.43,1.45,1.55,1.45\}$. Note that here, 45\% performance benefits for CHOLMOD in ``Arma Roller'' is equal to saving 7.86 hours which is achieved by simple integration of Parth without numerical side effect.  }

    \begin{tabular}{c|c c c c c c c}
    \midrule
     Parth + Tool & Step & Dolphin Funnel & Ball Mesh Roller & Mat On Board & Rods Twist & Squeeze out &  Arma Roller\\
     \hline
    \rowcolor{cyan!25}
     \cellcolor{white} & Symbolic time(s) & 3261 \ra 1424 & 2623 \ra 1134 & 1551 \ra 914 & 8635 \ra 3896 & 9098 \ra 4425 & 35357 \ra 15714\\
    \multirow{-2}{*}{ MKL} & Numeric time(s)& 932 \ra 923 & 1085 \ra 1075& 485 \ra 457 & 2800 \ra 2751  & 2543 \ra 2397 & 30703 \ra 32335     \\
    \midrule
    \midrule
        \rowcolor{cyan!25}
          \cellcolor{white}  & Symbolic time(s) & 1773 \ra 490 & 1465 \ra 369 & 985 \ra 483 & 5363 \ra 2481 & 5540 \ra 2119 & 25331 \ra 6913\\
         \multirow{-2}{*}{Accelerate} & Numeric time(s)& 706 \ra 709 & 861 \ra 861 & 412 \ra 412 & 2578 \ra 2602 & 2231 \ra 2251 & 51871 \ra 51098\\
    \midrule
    \midrule
            \rowcolor{cyan!25}
              \cellcolor{white}  & Symbolic time(s) & 2857 \ra 597 & 2353 \ra 453 & 1600 \ra 558 & 8552 \ra 2764  & 9012 \ra 2422 & 35021 \ra 6696\\
              \multirow{-2}{*}{CHOLMOD} & Numeric time(s) & 4200 \ra 4106 & 4511 \ra 4432 & 1663 \ra 1727 & 9098 \ra 9385 & 9483 \ra 9499 &  58282 \ra 57492\\
    \midrule
    \end{tabular}

    \label{tab:Evaluation:ParthEffectedRuntime}
\end{table*}

The key takeaway is the runtime savings achieved by integrating Parth. For example, using CHOLMOD for the "Arma Roller" simulation on an Intel processor offers numerical stability. With a solve speedup of 1.45x, the 45\% reduction in runtime results in a savings of 7.86 hours. Given the minimal overhead required to integrate Parth (three lines of code) and its consistent numerical stability, many practitioners can easily gain these performance benefits without needing to refactor an entire computational pipeline, which is often not straightforward. This reiterates Parth's objective to provide high-performance Cholesky solvers for a large group of practitioners instead of offering highly optimized, domain-specific solvers for specific applications.

Finally, we want to emphasize that, as computational pipelines—such as the one shown in IPC—are complex, accelerating the entire pipeline requires a deep understanding of its components and significant engineering effort. For example, in IPC simulations, it is well known that the CCD component is also a major bottleneck, and there have been substantial efforts to accelerate this component, such as the one presented in~\cite{IPCApprox2024}. However, as future work continues to improve the performance of these components, enhancements in the linear solver will become increasingly impactful.

\subsection{Remeshing: Benchmark}\label{sec:Remesher:Benchmark}


In this section we next evaluate how Parth performs in remeshing applications where sparsity patterns change due to localized updates in mesh geometry and topology which result in both adding and deleting of the nodes in G and added and removed edges. Here we test with a remeshing pipeline in which we select remeshing patch regions on triangle meshes, remesh with \citet{RemesherBotsch} to alter patch structure, and then apply global discrete Laplacian operator~\cite{libigl_laplace}. For a comprehensive analysis, remeshing operations are chosen to cover patches with a range of sizes comprising {1\%, 5\%, 10\%, 20\%} of the faces. Each patch is created around a randomly selected face ID. Here, per each mesh, the same face ID is used for comparison between different linear solvers. As a result, comparisons across linear solvers for a specific mesh and patch are consistent with identical computation. For each patch size, fifty such samples are chosen, covering different areas of a surface mesh, forming 200 sparse linear solves per mesh. After applying each patch remeshing, we reset the mesh and select another patch to measure Parth’s reuse capability for each patch size. This testing is applied across all the meshes in the \citet{odedstein-meshes} repository with the exclusion of the overly simple Cube mesh.

\begin{figure}
\centering
\includegraphics[width=\columnwidth]{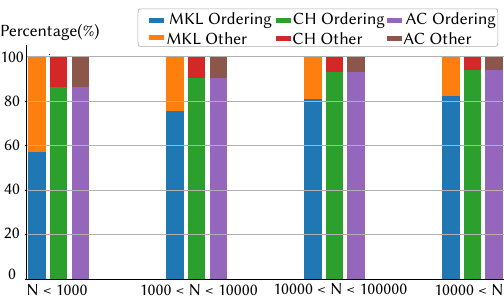}
\caption{{\bf Remeshing: Fill-reducing ordering is the bottleneck for all symbolic analyses.} Here we summarize a bottleneck analysis for the symbolic steps of all three Cholesky solvers across all 20 Meshes ranging from 642 DOFs to approximately 1.6m DOFs in our benchmark. For each simulation and corresponding tool, both the fill-reducing ordering time and all additional symbolic analysis time are recorded. Here, "CH", and "AC" respectively denote CHOLMOD and Accelerate LLT runtimes. "Other" categories summarize all other symbolic analysis costs per Cholesky solver with tasks that vary depending on the Cholesky solver method.}
\label{fig:Evaluation:Remesher_SymBottleneck}
\end{figure}

In addition to testing an application with changing matrix size, this benchmark additionally allows evaluation of Parth's performance on triangular meshes, and on linear solves with a problem dimension of $DIM=1$ (recalling for IPC, it is 3). As a result, in this evaluation, Parth directly employs the graph dual of the input matrix instead of the coarsened version where the graph nodes associated with a single DOF are combined. As a result, this benchmark does not require the application of Parth's compression.

\subsection{Remeshing: Bottleneck Analysis}

Figure~\ref{fig:Evaluation:Remesher_SymBottleneck} provides an analysis of the bottleneck of the symbolic stage. Here, we divide the meshes into separate groups based on their number of DOFs to show the effect of mesh sizes on the symbolic analysis stage. 
As we increase mesh size, and consequently increase the size of the graph dual of the Laplacian operator, the fill-reducing ordering performance overhead becomes more prominent. For example, for small meshes, fill-reducing ordering is 57\% of the symbolic analyzing, and for Large meshes, this overhead increases to 82\% for the MKL solver. As a result, for large-scale problems, optimizing symbolic analysis is clearly important for the Cholesky solve pipeline. Note that these results are consistent with our analysis of the IPC benchmark.

\begin{figure}[tp]
\centering
\includegraphics[width=\columnwidth]{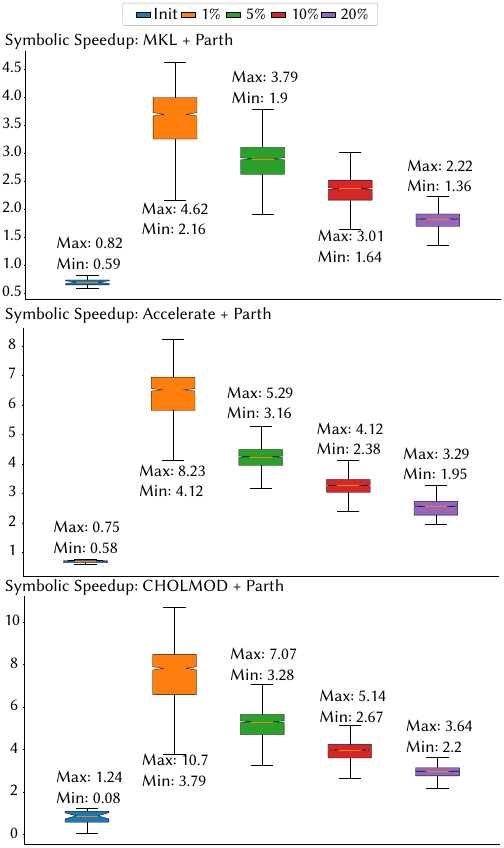}
\caption{{\bf Remeshing: Symbolic performance impact across Patch benchmark.} Here the performance impact of Parth on symbolic analysis of the Patch pipeline is shown. As it is shown, the initialization step is slowed down to 0.7x performance on average, due to Parth overhead. However, after that, due to the reuse a significant boost in performance is observed. Note that as expected, the reuse performance benefits are reduced due to the more aggressive changes in the mesh structure resulting from using the remesher.}
\label{fig:Evaluation:Remesher_SummaryWhisker_Sym}
\end{figure}

\subsection{Remeshing: Total Linear Solve Performance}

Here we discuss only the symbolic analysis and total per-solve linear solve performance of the remeshing benchmark. However, consistent with our IPC benchmark analysis, we also provide our numerical performance analysis of the remeshing benchmark in supplemental. For this set of analyses, we divide the performance data into five categories. The first category is the initialization step, while the next four categories reflect performance with respect to patch size. This categorization allows us to offer further insight into Parth's initialization cost without the compression phase, making the overhead of building the HGD data structure more visible. Furthermore, as shown in the teaser, a 2\% selection of the surface mesh is not considered small. By increasing these patch sizes to 20\%, we aim to provide insight into more challenging situations where the changes are more drastic, demonstrating Parth's performance in terms of both reuse capability and fill-reducing quality.


\begin{figure}[tp]
\centering
\includegraphics[width=\columnwidth]{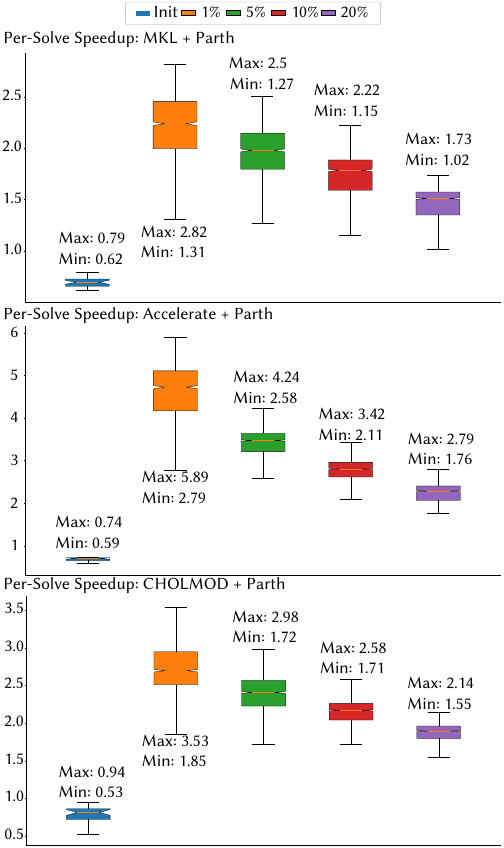}
\caption{{\bf Remeshing: Per linear solve speedup when using remesher.} Here, after each re-meshing, the per-solve speedup is shown. As we expected, when the numerical time comprises the smaller ratio of a solve runtime, Parth speedup is more significant. For example, since the Accelerate framework has a very fast numerical performance due to its fast processor, Parth shows up to 5.89x speedup with a median of 4.6x for a Patch size of 1\% of the total faces.}
\label{fig:Evaluation:Remesher_SummaryWhisker_PerSolve_Total}
\end{figure}

Figure~\ref{fig:Evaluation:Remesher_SummaryWhisker_Sym} shows Parth's effect on the symbolic analysis. As expected, since $DIM=1$ for this computational pipeline, Parth reduces performance in the initialization step. However, since the HGD computation is also reused in the fill-reducing ordering computation, the average performance speed is 70\% of the baseline fill-reducing routine. After the initialization step, Parth consistently provides speedup. For example, in CHOLMOD with patch sizes of 1\%, Parth achieves an order-of-magnitude speedup. It is important to note that as the size of the patches on the mesh increases, Parth’s performance decreases due to less temporal coherence and fewer opportunities for reuse, i.e., the changes become less gradual. In patch sizes of 20\%, for instance, we observe up to a 3.36x speedup, which is smaller than the order-of-magnitude speedup achieved with patch sizes of 1\% of the surface faces in symbolic analysis.

The performance benefits in symbolic analysis lead to an overall improvement in total Cholesky solve time, as shown in Figure~\ref{fig:Evaluation:Remesher_SummaryWhisker_PerSolve_Total}. As expected, this performance boost mirrors the trends seen in symbolic analysis, resulting in up to a 5.89x, 3.55x, and 2.82x speedup compared to Accelerate, CHOLMOD, and MKL, respectively.

\section{Limitations and Future Work}

Here, we discuss the key limitations of Parth that guide future research directions. Specifically, we discuss Parth's domain of applicability, its ability to maintain high-quality reuses under mesh deformation, and the challenges of alternative lagging approaches.

\subsection{Domain of Applicability}

Currently, Parth is designed to provide performance benefits by reusing symbolic analysis computations; however, its advantages are limited to scenarios where symbolic analysis—specifically fill-reducing ordering—is the primary bottleneck. When the computational bottleneck lies elsewhere, Parth’s improvements are not applicable. We plan to address this limitation in future work by adding adaptive numerical acceleration techniques to Cholesky solves.

\subsection{Widespread changes throughout the domain} 

Parth's core idea is to confine changes locally and enable reuse. Thus, if locality is not present, Parth fails to provide performance benefits. In the supplemental material, we provide a detailed analysis of Parth's reuse on the “Mat Twist” simulation from the IPC benchmark, where numerous self-collisions, spread throughout the domain, cause Parth to fail to provide reuse in some of the frames.

\subsection{Preservation of High-Quality Reuse}

Another key limitation worth investigating is how long Parth can preserve high-quality fill-reducing ordering performance, given that computing such orderings typically requires global information. To assess this, we apply remeshing to 1,000 patches (each comprising 1\% of mesh faces) on the triangular meshes evaluated in our Remeshing benchmark. Each patch is selected around face IDs not chosen in the previous patch using a uniform distribution, allowing us to sample different parts of the meshes uniformly. We then evaluate how many such 1\% patch remeshings can be applied before the numerical performance drops below 85\% of the baseline—indicating how often Parth can reuse high-quality ordering without a full recomputation. Our choice of a 15\% threshold is based on per-solve numerical speedup results (see supplemental), which show that a 15\% difference in numerical performance is typical for different high-quality orderings. For this test, the aggressive reuse algorithm (see supplemental) is activated to prevent full recomputation when a patch on the root separator changes. The results of this experiment provide insight into how much deformation from the original mesh Parth can tolerate before quality degrades.

\begin{figure}[tp] 
\centering 
\includegraphics[width=\columnwidth]{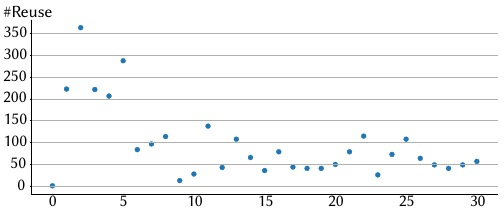} 
\caption{\textbf{Performance Degradation over Sequential Remeshings.} 
We evaluate how long Parth maintains high-quality fill-reducing ordering under continuous mesh modifications. For each mesh from~\cite{odedstein-meshes}, we perform a sequence of 1,000 remesh operations, each affecting 1\% of the mesh faces, followed by a discrete Laplacian smoothing. After each remeshing, we measure the numerical runtime of CHOLMOD factorization integrated with Parth ($t_p$) and compare it to the baseline runtime without Parth ($t$) by computing the relative difference $(t_p - t) / t$. The figure plots, for each mesh, the number of remeshings that can be applied before this relative difference exceeds 15\%, indicating that Parth’s performance falls below 85\% of the baseline due to degraded ordering quality. The results reveal that, on average, performance drops after 94 remesh operations, with a median threshold reached after approximately 64\% of the mesh has been modified.} 
\label{fig:Limitation:Limit}
\end{figure}

As illustrated in Figure~\ref{fig:Limitation:Limit}, the x-axis shows mesh IDs from~\cite{odedstein-meshes} (sorted by degrees of freedom), and the y-axis indicates the number of reuses. For some meshes, Parth maintains high-quality ordering even after 350 patches. On average, performance drops below 85\% after 94 distinct patch remeshings—corresponding to 94\% of the mesh being modified. However, results vary depending on mesh shape and deformation extent. The median result indicates that after 64.42\% of the mesh has changed, Parth needs to be reset to obtain new fill-reducing ordering information, as the underlying geometry has shifted significantly.

\subsection{Lagging Approaches as Alternatives}

Finally, we consider lagging approaches as an alternative to Parth. To the best of our knowledge, no lagging approach supports reusing the fill-reducing ordering when the number of matrix rows or columns changes. For cases where only the non-zero entries of the linear system change, one can reuse the full permutation vector via a heuristic that delays the order computation. In practice, these approaches are difficult to use because they require domain-specific knowledge and lack generality as the performance can degrade significantly if lagging approach is not appropriate. We provide a comprehensive analysis of two such lagging approaches in the supplemental, illustrating the challenges in designing them and offering insight into scenarios where they might replace Parth.
\section{Conclusion}
Our evaluation reveals that Cholesky solves face significant bottlenecks when the sparsity pattern is dynamic due to symbolic analysis overhead—a factor often ignored in the recent development of high-performance Cholesky solvers. Our study further pinpoints the fill-reducing ordering algorithm as a key contributor. This insight offers a clear path to mitigate the bottleneck while continuing to use state-of-the-art Cholesky solvers. Our extensive evaluation confirms that the Parth module enables successful, simple, and direct integration into all three state-of-the-art Cholesky solver libraries—without the need for any tuning and using only the input matrix (and a mapping when the matrix rows/columns change). The generality and ease of use of Parth align with practitioners’ objectives on using Cholesky solves, whether they prioritize accuracy and stability or face constraints on developing high-performance specialized iterative linear solvers. Building on this foundational capability, we next conclude how Parth improves performance by reducing the bottleneck in terms of both runtime and quality."

Our evaluation shows that integrating Parth into high-performance Cholesky solvers yields significant speedups and reliability improvements across various challenging linear benchmarks. By reusing computations across calls, Parth achieves up to a 255× speedup in fill-reducing ordering computation while preserving numerical stability and performance in the downstream Cholesky factorization and subsequent solve processes for all three high-performance direct solvers—namely, Intel MKL Pardiso~\cite{MKLPardiso}, CHOLMOD~\cite{Cholmod}, and the newly introduced Apple Accelerate sparse solver~\cite{Apple}. This integration leads to up to a 5.89× speedup in overall Cholesky solve time, outperforming the best speedups obtained by recent architecture-customized Cholesky solver updates~\cite{Apple}. As Cholesky solvers remain one of the core time-intensive computational routines in applications with dynamically changing sparsity, these performance gains represent a significant improvement enabled by a simple three-line module addition at the solver’s API level. These achievements not only address current challenges but also allow for many future enhancements that will further extend Parth’s capabilities.

Finally, to further extend Parth’s advantages, we plan to reuse additional parts of the symbolic analysis and develop faster numerical computation for Cholesky solves, leveraging the low-overhead symbolic analysis already achieved by Parth. Alongside these technical advances, our open-sourced Parth codebase invites community collaboration, enabling researchers and practitioners to build upon and refine our methods. We believe that Parth will play a significant role in the future advancement of high-performance Cholesky solvers for handling dynamically changing sparse linear systems efficiently, setting the stage for continued innovation, broader adoption, and collaborative development.

\begin{acks}
This work was supported by NSERC Discovery Grants (RGPIN‑2019‑06516 and RGPIN‑2023‑05120),  
the Canada Research Chairs Program,  
the Ontario Early Researcher Award,  
and the Digital Research Alliance of Canada (\url{www.alliancecan.ca}).  
\end{acks}

\bibliographystyle{ACM-Reference-Format}
\bibliography{references}



\end{document}